\documentclass[10pt]{amsart}
\usepackage{amsmath}
\usepackage{amsthm}
\usepackage{latexsym}
\usepackage{amssymb}

\newtheorem{theorem}{Theorem}
\newtheorem{proposition}{Proposition}
\newtheorem{lemma}{Lemma}

\newtheorem{definition}{Definition}

\newtheorem{remark}{Remark}
\newcounter{obsctr}

\renewcommand{\theequation}{\thesection.\arabic{equation}}

\begin{document}
\def\A {{\mathcal{A}}}
\def\D {{\mathcal{D}}}
\def\R {{\mathbb{R}}}
\def\N {{\mathbb{N}}}
\def\C {{\mathbb{C}}}
\def\Z {{\mathbb{Z}}}
\def\l {\ell}
\def\tm {\tilde{m}}
\def\ml {multline}
\def\multiline {\multline}
\def\lessim {\lesssim}
\def\phi{\varphi}
\def\epsilon{\varepsilon}
\def\olm{\overline{L_m}}
\def\ol{\overline{L}}
\def\oz{\overline{z}}
\def\be{\begin{equation}}
\def\\[{\begin{equation}}
\def\ee{\end{equation}}
\def\\]{\end{equation}}
\title{Singular sums of
squares of degenerate vector fields}    
\author{Antonio Bove}
\address{Dipartimenot di Matematica, Universit\`a
di Bologna, Piazza di Porta San Donato 5, Bologna
ITALY}
\email{bove@unibo.edu}
\author{Makhlouf Derridj}
\address{5 Rue de la Juvini\`ere, 78350 Les Loges en
Josas, FRANCE}
\email{derridj@club-internet.fr}
\author{Joseph J. Kohn}
\address{Department of Mathematics, Fine
Hall, Princeton University, Princeton NJ 08540, USA}
\email{kohn@math.princeton.edu}
\author{David S. Tartakoff}
\address{Department of Mathematics, University
of Illinois at Chicago, m/c 249, 851 S.
Morgan St., Chicago IL  60607, USA}
\email{dst@uic.edu}
\date{\today}
\begin{abstract} In \cite{K2005}, J. J. Kohn 
proved $C^\infty$ hypoellipticity with loss of
$k-1$ derivatives in Sobolev norms (and at least that
loss in $L^\infty$) for the highly non-subelliptic
singular sum of squares 
$$ P_k=L\overline{L} +
\overline{L}|z|^{2k}L =-\overline{L}^*\overline{L} -
(\overline{z}^kL)^*\overline{z}^kL\hbox{ \quad with \quad}
L={\partial
\over
\partial z} + i\overline{z}{\partial \over
\partial t}.$$

In this paper, we prove hypoellipticity with loss
of $\frac{k-1}{m}$ derivatives in Sobolev norms for the
operator 
\begin{equation} P^F_{m,k}=L^F_m\overline{L^F_m} +
\overline{L^F_m}\,|z|^{2k}L^F_m \hbox{ \; with \;}
L^F_m={\partial
\over
\partial z} + iF_{z}{\partial \over
\partial t},\end{equation}
with $F(z,\oz)$ such that 
\begin{equation}F_{z\oz} = |z|^{2(m-1)}g,\;
g(0)>0,\hbox{ so that }  F_z =
\oz|z|^{2(m-1)}h\end{equation}  whose
prototype, when
$mF(z,\oz)=|z|^{2m},$ is 
\begin{equation} P_{m,k}=L_m\overline{L_m} +
\overline{L_m}\,|z|^{2k}L_m, \;\;
L_m={\partial
\over
\partial z} + i\overline{z}|z|^{2(m-1)}{\partial \over
\partial t},\end{equation}
for which the underlying manifold is of finite
type. 

We give two proofs: the first using a fairly
rapid derivation of an {\it a priori} estimate
analogous to that used by Kohn 
in \cite{K2005}:
\begin{equation}\|\phi u\|_0 \leq
C\|\tilde{\phi}P_{m,k}^F v\|_\frac{k-1}{m} +
C\|u\|_{-\infty}
\end{equation}
(for all $u \in C_0^\infty$ with $\tilde\phi \equiv 1$
near supp $\phi$), after deriving this estimate in the
first part of the paper; 
 the second uses the far more rapidly
derived estimate of
\cite{Ta2005} and
\cite{DT2005} (where analytic hypoellipticity for
$P_k$ and
$P_{m,k}$ are also proved): $\forall v \in C_0^\infty$ of small
support,
\begin{equation}\|v\|_{-\frac{k-1}{2m}}^2 + \|\ol
v\|_0^2+\|\oz^k Lv\|^2_0\leq C|(P_{m,k}^Fv,v)_{L^2}| +
C\|v\|_{-N}^2.\end{equation}
We also prove, along the
way, analytic hypoellipticity for $P^F_{m,k}.$ For
\begin{equation}F(z,\oz)=f(|z|^2),
\end{equation}
we show that these estimates are optimal.
\end{abstract}
\maketitle
\pagestyle{myheadings}
\markboth{Bove, Derridj, Kohn \& Tartakoff}
{Singular Sums of Squares of Vector Fields}
\section{Introduction and statement of theorems}
\renewcommand{\theequation}{\thesection.\arabic{equation}}
\setcounter{equation}{0}
\setcounter{theorem}{0}
\setcounter{proposition}{0}  
\setcounter{lemma}{0}
\setcounter{corollary}{0} 
\setcounter{definition}{0}
\setcounter{remark}{0}

In his recent paper, \cite{K2005}, J. J. Kohn
exhibited a sum of squares of complex vector fields
which satisfied the bracket condition but which was
not subelliptic; nonetheless, he showed that the
operator was hypoelliptic, though with a large loss of
derivatives. His example was:
$$ P_k=L\overline{L} +
\overline{L}|z|^{2k}L =-\overline{L}^*\overline{L} -
(\overline{z}^kL)^*\overline{z}^kL\hbox{ \quad with \quad}
L={\partial
\over
\partial z} + i\overline{z}{\partial \over
\partial t}.$$

The {\it a priori} estimate Kohn established is a strong
one and in this case (since the operator is
independent of the variable $t,$) leads virtually at
once to the hypoellipticity of $P_k$: for any $s,$ there
exists a constant $C_s$ such that for all smooth $u$ and
any pair of cut-off functions $\phi, \tilde{\phi}$ with
$\tilde{\phi}\equiv 1$ near supp $\phi,$
\begin{equation}\|\phi u\|_s \leq
C_s\|\tilde{\phi}P_{k} u\|_{s+k-1} +
C_s\|u\|_{-\infty}
\end{equation}
Here the last norm stands for a norm of arbitrarily low
order, with the constant preceeding it possibly
depending on the order of that norm, and $u$ assumed to
be of (possibly large) compact support.

Subsequently, in
\cite{DT2005}, M. Derridj and D. S. Tartakoff
proved analytic hypoellipticity for $P_k$
using rather different methods, namely they established
an inequality for functions $v$ of small support, hence
an estimate which did not require explicit cut-off
functions, reserving the necessity of localizing an
actual solution to a neighborhood of a point to the
proof of (analytic) hypoellipticity: for any $s,$ there
exists a constant $C_s$ such that for all $v\in
C_0^\infty$ of small support, 
\begin{equation}
\|v\|_{s-\frac{k-1}{2}}^2+\|\ol v\|_s^2 + \|\oz^k L
v\|_s^2\leq C_s |(P_k v, v)_s| + C_s\|v\|_{-\infty}
\end{equation}
which of course yields the previous estimate at once
without the cut-off functions but only for $u$ already
known to have (small) compact support. 

This paper was partly motivated by the effort to
understand the relationship between these
estimates, partly to obtain a simpler (or at least more
concise) derivation of the former, and finally to
generalize these results where possible. 

In \cite{Ta2005}, the fourth author had already
sharpened the methods of \cite{DT2005} to include the
example of the operator 
\begin{equation} 
P_{m,k}=L_m\overline{L_m} +
\overline{L_m}\,|z|^{2k}L_m \hbox{ \; with \;}
L_m={\partial
\over
\partial z} + i\overline{z}|z|^{2(m-1)}{\partial \over
\partial t}
\end{equation}
based on the tangential vector fields to a domain in
$\C^2$ of finite type; the technical work was
heavily dependent on the methods of \cite{DT1988}. 

Both \cite{Ta2005} and \cite{DT2005} include proofs of
$C^\infty$ - hypoellipticity by `truncating' the proofs
of analytic hypoellipticity, hence use 
the entire
machinery that has come to be known as $(T^p)_\phi$
since \cite{Ta1978}. 

In this paper, we consider the more general operator 
\begin{equation} P^F_{m,k}=L^F_m\overline{L^F_m} +
\overline{L^F_m}\,|z|^{2k}L^F_m \hbox{ \; with \;}
L^F_m={\partial
\over
\partial z} + iF_{z}{\partial \over
\partial t},\end{equation}
with 
\begin{equation}\label{F} F_{z\oz} = |z|^{2(m-1)}g,\;
g(0)\neq 0,\; F_z =
\oz|z|^{2(m-1)}h\end{equation}
 whose prototype, when
$F(z,\oz)=|z|^{2m}/m,$ is the operator $P_{m,k}$
discussed above.

Here we establish two families of 
estimates for 
$P^F_{m,k},$ and prove the optimality of these
estimates under the additional restriction
\begin{equation}\label{F'}F(z,\oz) =
f(|z|^2).\end{equation} 

We will then use one of the estimates
 to prove
$C^\infty$ hypoellipticity with precise loss of
$\frac{k-1}{m}$ derivatives and the other to prove 
$C^\omega$ hypoellipticity and to give another proof of
$C^\infty$ hypoellipticity with the prescribed loss.

A note on the norms used is in order. All of our norms
and derivations are done in $L^2(z,\oz)\times H^s(t).$ 
There are several reasons for this. First, Proposition
1.1 could, for $s=0,$ trivially have the norm on the
left replaced with the full $-\frac{k-1}{2m}$ norm,
then as mentioned below, using a cut-off in $\tau$ dual
to $t$ which tends to the identity, one can prove
easily that since $\partial_t$ commutes with $P,$ high
$t$ derivatives of the solution belong to
$H^{s-\frac{m-1}{m}}$ (in $t$) provided this is true
of $Pu$ in $H^s$ norm. 

But the whole classical theory of pseudo-differential
operators and wave front sets allows us to
microlocalize the consideration of
hypoellipticity. For it is clear that if
$z\neq 0,$ the operator is elliptic and
hence even analytic hypoelliptic, and
gains two derivatives. For $z$ close to
zero, one must look in the cotangent
space, $(z,t;\zeta, \tau)$ which, in the complement of
$(z,t;0,0)$ we write as the union of overlapping cones:
the cones
$\Gamma^\pm$ contain
$\tau = + 1,
\zeta = 0$ and $\tau = - 1, \zeta =
0$ respectively, while the ``elliptic" cone $\Gamma^0$
contains $\tau = 0.$ In $\Gamma^0,$ the operator $P$ is
also elliptic, since this is true of 
$L\ol.$ In $\Gamma^-,$ the operator $P$ is maximally 
hypoelliptic and hence is
subelliptic with loss of
$1/2m$ derivatives ($\|Lv\|$ is bounded by $\|\ol v\|$
there, hence the operator is maximally hypoelliptic,
which means that the real and imaginary parts of $L$
and $\ol$ are bounded by $P,$ and by H\"ormander's
condition, subelliptic and one has the estimate of
Proposition 1.1 with the $H^{s+\frac{1}{2m}}$ norm
on the left), hence is (microlocally) hypoelliptic
with a gain of $\frac{1}{2m}$ derivatives in that region
by conventional arguments.  

It is only in the positive
cone that all of this work is necessary, and there in
addition to having estimates such as Lemma 2.3 below,
we also know (as we would in $\Gamma^-$ as well) that 
$|\zeta|\leq C|\tau|$ so that estimating high
derivatives in $t$ will yield control in all
directions. 

In
the two Propositions which follow, the notation
$A\lesssim B$ will mean that $A\leq CB$ with $C$
uniform in $v\in C_0^\infty$ and locally so in $s,$
and $F$ is assumed to satisfy the conditions of (2.1)
above.

\begin{proposition} For $v$ of small support, 
\begin{equation}\label{ipape}\|v\|^2_{s-\frac{k-1}{2m}}
+
\|\overline {L_m^F} v\|_s^2 +
\|\oz^kL_m^Fv\|_s^2 \lesssim 
|(P_{m,k}^Fv,v)_s|+\|v\|_{-\infty}^2,
\end{equation}
\end{proposition}
\begin{proposition}
For any pair of cut-off functions 
$\phi, \tilde{\phi}$ with
$\tilde{\phi}\equiv 1$ near supp $\phi,$ and for $u$ of
support in a fixed (not necessarily small) compact set,
\begin{equation}\label{coape}\|\phi u\|_s^2 \lesssim
\|\tilde\phi P_{m,k}^Fu\|_{s+\frac{k-1}{m}}^2 +
\|u\|_{-\infty}^2
\end{equation}
\end{proposition}
\begin{proposition} For the case
$F(z,\oz) =
f(|z|^2)b(z,\oz), b(0)\neq 0$ the loss in
Propositions 1.1 and 1.2 cannot be improved.
\end{proposition}
\begin{theorem} $P_{m,k}^F$ is locally 
hypoelliptic with loss of $\frac{k-1}{m}$ derivatives:
$P_{m,k}^Fu\in H^s \implies u\in
H^{s-\frac{k-1}{m}}$.
\end{theorem}

In the sequel, we will write
$L$ for $L_m^F,\ol$ for $\overline{L_m^F},$ and $P$ for
$P_{m,k}^F.$

\section{Preliminary Observations and Lemmas}
\renewcommand{\theequation}{\thesection.\arabic{equation}}
\setcounter{equation}{0}
\setcounter{theorem}{0}
\setcounter{proposition}{0}  
\setcounter{lemma}{0}
\setcounter{corollary}{0} 
\setcounter{definition}{0}
\setcounter{remark}{0}

The first observation concerns the apparent difference
between the two {\it a priori} estimates in the
two Propositions above and their use. The second
estimate explicitly introduces a second cut-off
function, although, as we shall see, except for the
last term, the function $\tilde\phi$ may be replaced by
certain derivatives of $\phi.$ That is, except for a
norm of sufficiently low order, we may control the
terms on the right by derivatives of the given
localizing function. In fact, {\it the same is true in
the proof of (analytic) hypoellipticity using the first
estimate} - we proceed with a balanced
localization $(T^p)_\phi$ of high derivatives in
$T=\partial_t$ and encounter errors expressed as
derivatives of the localizing function we start with
and then at a certain point (in this case a fraction of
the derivatives we seek to estimate), we are forced to
introduce a cut-off function with strictly larger
support and to construct a whole new balanced sum 
$(T^{\tilde{p}})_{\tilde\phi}$ around this new
localizing function - and for the analyticity proof we
need to control these supports in a very precise way.

It
is not at all clear how to pass from one setting to the
other - neither estimate trivially implies the other
and the proofs of hypoellipticity are not trivially
comparable, but they do seem to contain the same
elements.

Our first technical
observation concerns the dependence of localizing
functions on
$z.$ In order to localize to a neighborhood of
$0,$ we may take a product of a function of
$z,\oz$ of small support but identically equal to
one near the origin in
$\C$ with another function of $t$ only, again taken to
be of small support. Whenever the first of these
functions is differentiated, the resulting function is
supported away from $z=0,$ hence in a region where the
operator $P$ is in fact subelliptic and hence far
better behaved. We shall ignore such
regions and thus take all localizations to be
functions of $t$ only. 

To make the proofs of Propositions 1.1 and 1.2 flow
more smoothly, we prepare some easy lemmas which will be
used repeatedly in the sequel. By integration by parts
and shifting powers of $z$ from one side of an inner
product to the other, these lemmas, especially
Lemma 2.4 and 2.5, which are often used, express the
obvious fact that by grouping one power of $z$ and a
fractional power of $\Lambda_t,$ effectively a
fractional power of $\partial_t,$ as a unit, say
$A=z\Lambda_t^\rho,$ one may move powers of $A$ 
from one side of an inner product to the other. In
all of these lemmas,
$w$ will denote a smooth function of (small) compact
support and the superscript `+' will indicate that the
function has been microlocalized to the positive cone
for the symbol of $\partial_t.$ The estimates are
locally uniform in
$s.$ 
\begin{lemma} $\|Lw\|_{s-1/2}\lesssim \|\ol w\|_{s-1/2}
+ \|z^{m-1}w\|_s.$
\end{lemma}
\begin{proof} Integration by parts since $[L,
\ol] = -2iF_{z\oz}\partial_t, |F_{z\oz}|\lesssim
|z|^{2(m-1)}.$
\end{proof}
\begin{lemma} $\|z^{m-1}w\|_{s+1/2}
\lesssim \|Lw\|_{s}+
\|\ol w\|_{s}.$
\end{lemma}
\begin{proof} Integration by parts since $[L,
\ol] = -2iF_{z\oz}\partial_t$ and $F_{z\oz}\geq
c|z|^{2(m-1)}.$
\end{proof}
\begin{lemma} $\|\ol w^+\|_s +
\|z^{m-1}w^+\|_{s+1/2}\lesssim
\|Lw^+\|_s + \|w\|_s.$
\end{lemma}
\begin{proof} The same identity where the symbol of
$-2i\partial_t$ has the appropriate sign.
\end{proof}
\begin{lemma}$\|z^rw\|_\mu \leq
s.c.\|z^{r-n_1}w\|_{\mu-n_1\rho}
+l.c.\|z^{r+n_2}w\|_{\mu+n_2\rho}, \;n_1\leq r.$
\end{lemma}
\begin{proof} Let $A=z\Lambda_t^\rho.$ Then for example
$$\|A^{r}w\|_\mu^2 = (A^{r-n_1}w, A^{r+n_1})_\mu \leq
s.c.\|A^{r-n_1}w\|^2_\mu+l.c.\|A^{r+n_1}w\|^2_\mu$$ but
then the second of these terms may be related to lower
and higher powers of $A,$ and the result follows.
\end{proof}
\begin{lemma}$\|z^rw\|_\mu \leq
l.c.\|z^{r-n_1}w\|_{\mu-n_1\rho}
+s.c.\|z^{r+n_2}w\|_{\mu+n_2\rho}, \;n_1\leq r.$
\end{lemma}
\begin{proof} Completely analogous.
\end{proof}
\begin{lemma}$\|w\|\lesssim \|z\ol w\|+\|zLw\|.$
\end{lemma}
\begin{proof} This is the subelliptic multiplier
argument: $$\|w\|^2 = |([L,z]w,w)| \leq |(w,\oz\ol w)| +
|(zLw,w)| $$
$$\leq s.c.\|w\|^2 + l.c.(\|z\ol w\|^2 +
\|zLw\|^2).$$
\end{proof}
\begin{lemma}$\|\phi u\|_0\lesssim \|z\ol \phi u\|_0 +
\|\oz L\phi u\|_0.$
\end{lemma}
\begin{proof} This is the previous lemma with $w=\phi
u.$
\end{proof}
\begin{lemma}$\|\phi u\|_0\lesssim \|z\phi \ol u\|_0 +
\|\oz \phi L u\|_0 + \|z^{2m}\phi'u\|_0.$
\end{lemma}
\begin{proof} This is just the observation that
$|[L,\phi]u|\sim |z|^{2m-1}|\phi_t|.$
\end{proof}
\begin{lemma}$ \|z\phi' u\|_0 \lesssim
\|
 \ol z\phi'u\|_{-\frac{1}{2m}} +
\|L z \phi'u\|_{-\frac{1}{2m}} +
\|z \phi' u\|_{-\frac{1}{2m}}$
\end{lemma}
\begin{proof} This is just the observation that the
vector fields $L$ and $\ol,$ or rather their real and
imaginary parts, satisfy the (real) bracket condition
and hence form a subelliptic system in the usual sense
with $\epsilon = 1/2m,$ and then the whole subelliptic
estimate is lowered by $1/2m.$
\end{proof}
\section{Proof of Proposition 1.1}
\renewcommand{\theequation}{\thesection.\arabic{equation}}
\setcounter{equation}{0}
\setcounter{theorem}{0}
\setcounter{proposition}{0}  
\setcounter{lemma}{0}
\setcounter{corollary}{0} 
\setcounter{definition}{0}
\setcounter{remark}{0}

To prove Proposition 1.1, the {\it a priori} estimate
on compactly supported functions, we set
$r=-{k-1\over 2m}$ and
$\tau = {1\over 2m}.$ Note that $r$ need not be
negative, but
$r-\tau = -{k\over 2m}\leq 0.$ Then we have,
since $r\leq \tau$: 
$$\|v\|_r^2 = ((Lz)\Lambda_t^rv, \Lambda_t^rv)  = 
(z\Lambda_t^{2r}v, \overline{L}_mv) -
(L\Lambda_t^{r-\tau}v,
\overline{z}\Lambda_t^{r+\tau}v)$$
$$\leq 
C\{\|\ol v\|^2+
s.c.\|L\Lambda_t^{r-\tau}v\|^2 +
l.c.\|z\Lambda_t^{\tau} (\Lambda_t^r v)\|^2\}$$
$$\leq C\{\|\ol v\|^2+
l.c.\|z\Lambda_t^{\tau} (\Lambda_t^r v)\|^2 + 
s.c.\|z^{m-1}\Lambda_t^{{1\over 2}-\tau}(\Lambda_t^r
v)\|^2\}$$
$$\leq C\{\|\ol v\|^2+
l.c.\|z\Lambda_t^{\tau} (\Lambda_t^r v)\|^2 + 
s.c.\|z^{m-1}\Lambda_t^{(m-1)\tau}(\Lambda_t^r
v)\|^2\}.$$

When $m=1,$ this last term is just
$s.c. \,\|v\|_r^2$ but for $0<a\leq m-1+k,$ we use
Lemma 2.4 in the form
$$\|z^a\Lambda_t^{a{\tau}}w\|^2 \leq s.c. \|w\|^2 +
l.c.
\|z^{m-1+k}\Lambda_t^{(m-1+k){\tau}}w\|^2$$
$$= s.c. \|w\|^2 + l.c.
\|z^{m-1+k}\Lambda_t^{1\over
2}(\Lambda_t^{(m-1+k){\tau}-{1\over 2}}w)\|^2$$ 
with $w=\Lambda_t^rv$ twice, once for $a=1$ and once
for $a= m-1.$  

Inserting this in the estimate above for
$\|v\|^2_r,$ we find 
$$\|v\|_r^2 \leq C\{ \|\ol v\|^2 +
\|z^{m-1+k}\Lambda_t^{1\over 2}\Lambda_t^{(m-1+k)\tau
-{1\over 2}+r} v
\|^2\}$$
$$= C\{ \|\ol v\|^2 +
\|z^{m-1+k}\Lambda_t^{1\over 2}
v\|^2\}.$$
since $(m-1+k)\tau - {1\over 2}+r=0.$

On the other hand, we have by Lemma 2.2, 
$$\|z^{m-1+k}\Lambda_t^{1\over 2}v\|^2 \leq
C\{\|\ol v\|^2
+\|\overline{z}^kLv\|\},$$ (which one proves
from the Lemma with the additional term
$\|z^{k-1}v\|^2$ on the right and then,
writing $z^{k-1}\sim [L, z^k],$ absorbs this
term by the other two). Thus we arrive at 
$$\|v\|_r^2 \leq C\{ \|\ol v\|^2
+\|\overline{z}^kLv\|^2\}=C|(Pv,v)|\leq
l.c.\|Pv\|_{-r}^2+s.c.\|v\|_r^2$$
or 
$$\|v\|_{-\frac{k-1}{ 2m}}^2+ \|\ol v\|^2
+\|\overline{z}^kLv\|^2\lesssim
C\|Pv\|^2_{\frac{k-1}{2m}}, v\in C_0^\infty, \qed$$

\section{Proof of Proposition 1.2. The Case $k=1$}
\renewcommand{\theequation}{\thesection.\arabic{equation}}
\setcounter{equation}{0}
\setcounter{theorem}{0}
\setcounter{proposition}{0}  
\setcounter{lemma}{0}
\setcounter{corollary}{0} 
\setcounter{definition}{0}
\setcounter{remark}{0}

For $k=1$ we will establish the estimate (for $u$ of
small support near $z=0$) using only Lemmas 2.3 and
2.8:
\begin{multline}\label{k=1}\|\phi u\|_s^2 + \|\phi\ol
u\|_s^2 + \|\oz \phi Lu\|_s^2\equiv \sum_{j=1}^{3}
(LHS)_j
\lesssim \\
\lesssim \sum_{j=0}^N\|\phi^{(j)}P
u\|^2_{s-j/2} +
\|\phi^{(N)} u\|^2_{s-{N/2}} + \|u\|_{-\infty}^2.
\end{multline}

Here and elsewhere, we will find the following
definition useful: 

\begin{definition} The designation ``RJ" (for
``Relative Junk") will apply to any multiple of any of
the terms $(LHS)_j$ that we are in the process of
estimating but with lower Sobolev index and possibly a
derivative on the localizing function - in other words,
to a term which will be iteratively estimated at the
end.
\end{definition}

For any value of $k,$ from Lemma 2.8, 
$$(LHS)_1\equiv\|\phi u\|^2 \lesssim \|z\phi \ol u\|^2
+
\|\oz
\phi Lu\|^2 + \|z^{m}\phi' u\|^2$$
(Lemma 2.8 even gives $z^{2m}$ in place of $z^m$). We
claim that this last term is RJ. To see this, Lemma 2.3 
tells us that
\begin{equation}\label{RJ1}\|z^m\phi'u\|^2_0 =
\|z^{m-1}\oz\phi'u\|_0^2
\lesssim
\|\oz L\phi'u\|^2_{-1/2} +
\|z\phi'u\|^2_{-1/2}
= RJ \end{equation}
(provided, as we will show, that we can estimate $\|\phi
u\|_0$ and
$\|\oz L\phi u\|_0^2$). Actually, in the next section
we will see even that $\|z\phi'u\|_0 \in RJ.$

So we have, modulo $RJ$ 
$$(LHS)_1\lesssim (LHS)_2 + (LHS)_3 \equiv \|\oz \phi L
u\|_0^2+\| \phi \ol u\|^2_0$$
$$=|-(\ol \phi^2 |z|^2 Lu,u) - (L\phi^2\ol u,u)| \lesssim
|(\phi^2 P u,u)|$$
$$ + |(F_z\phi\phi'\ol u,u)| +
|(F_{\oz}\phi\phi'|z|^2 Lu,u)| + |(\phi^2\oz Lu,u)|.  $$
These last two terms are easy to
handle: 
$$  |(F_z\phi\phi'\ol u,u)| +
|(F_{\oz}\phi\phi'Lu,u)| $$
$$
\lesssim s.c.\|\phi\ol u\|^2 + l.c.\|z^{2m-1}\phi' u\|^2 +
s.c.\|\oz L u\|^2$$
which are absorbed modulo the term $\|z^{2m-1}\phi' u\|$
which is RJ since $2m-1\geq m.$ 

Thus in all, in the positive cone, with $\tilde\phi \equiv
1$ near the support of $\phi,$ and for any $N,$
$$\|\phi u\|_s^2 + \|\phi\ol
u\|_s^2 +
\|\oz
\phi Lu\|_s^2
\lesssim \sum_{j=0}^N\|\phi^{(j)}P
u\|^2_{s-j/2} +
\|\phi^{(N)} u\|^2_{s-{N/2}}+\|u\|_{-\infty}^2.$$
Or, with $\tilde\phi \equiv 1$ near the support of $\phi,$ 
$$\|\phi u\|_s^2 
\lesssim \|\tilde\phi P u\|_s^2 +
\| u\|^2_{-\infty}.$$

\section{Proof of Proposition 1.2. The Case $k>1$}
\renewcommand{\theequation}{\thesection.\arabic{equation}}
\setcounter{equation}{0}
\setcounter{theorem}{0}
\setcounter{proposition}{0}  
\setcounter{lemma}{0}
\setcounter{corollary}{0} 
\setcounter{definition}{0}
\setcounter{remark}{0}

To prove Proposition 1.2 when $k>1,$ which is harder, we
cannot just use Proposition 1.1 with a cutoff function
$\phi$ in front of $v$ and then express the right hand
side in terms of $(\phi P v, \phi v)$ modulo acceptable
errors, since bracket of $P$ with $\phi$
introduces errors easily absorbed only when the
basic estimate is subelliptic, which here means
$k=0,$ the well-known case, or at least, in Kohn's
terminology, `no loss, no gain', namely the case
$k=1$ which we just considered.  

Instead, we proceed as follows. We will
establish again the class of ``Relative Junk Terms'',
denoted
$RJ,$ which are of the same form as those terms being
estimated but of lower Sobolev degree, and the
localizing function(s) may have received a derivative.
These will be treated recursively at the end, in a very
simple manner, but to see that a term is RJ one may
have to compare it to all eight terms below. 

The terms we want to estimate are eight in number, and
will be referred to as $(LHS)_j, j = 1, \ldots 8.$ In
estimating some the others will occur, generally with
small constants, but we set up a generic sum with
unknown coefficients $\sum_{j=1}^{j=8} C_j(LHS)_j.$
Specifically, we will establish, for suitable $C_j$ to
be determined relative to one another, 

$$C_1\|\phi u\| + C_2\|\oz \phi Lu\|+ C_3\|\phi \ol
u\|_{\frac{k-1}{2m}}+ C_4\|\phi\oz^{k} L
u\|_{\frac{k-1}{2m}}+ C_5\|z^{2k+m-2}\phi
u\|_{\frac{k-1}{m}}$$
\begin{equation}\tag{*}+C_6\|\phi
z^{2k-1}Lu\|_{\frac{k-1}{m}-\frac{1}{2}} +
C_7\|z^{m-1}\phi\ol u\|_{\frac{k-1}{m}}+C_8\|L \phi\ol
u\|_{\frac{k-1}{m}-\frac{1}{2}}\end{equation}
$$ =E=\sum_1^8 C_j(LHS)_j\leq C_9 \|\phi
Pu\|_{\frac{k-1}{m}} +
RJ .$$

In proving $(*)$ we will encounter errors from
microlocalization, errors which are supported in
regions where the regularity is well understood. As
these will be included in $RJ$ in any case, we will
omit explicit mention of terms of the form
$\|u\|_{-\infty}.$

\subsection{Estimating $(LHS)_1$}

Using Lemma 2.8 and then Lemmas 2.9
and 2.4:
$$(LHS)_1 \equiv \|\phi u\|_0 \lesssim
\|\phi\overline{L} u\|_0+
\|\overline{z}\phi L u\|_0 + \|z\phi' u\|_0$$
$$\lesssim
\|\phi\overline{L} u\|_0+
l.c.\|\oz^k \phi Lu\|_{\frac{k-1}{2m}} + s.c.\|\phi
Lu\|_{-\frac{1}{2m}} + \|z\phi' u\|_0.$$ 
For third term we write, using Lemma 2.1:
$$\|\phi Lu\|_{-\frac{1}{2m}} \lesssim \|z^{m-1}\phi
u\|_{-\frac{1}{2m}+\frac{1}{2}} +\|\phi
\ol u\|_{-\frac{1}{2m}}+\|\phi' u\|_{-\frac{1}{2m}}$$
(the last two terms are RJ) and from Lemma 2.4, 
$$ \|z^{m-1}\phi
u\|_{-\frac{1}{2m}+\frac{1}{2}}\lesssim s.c.\|\phi u\|_0 +
l.c.\|z^{m-1}\oz^k \phi
u\|_{\frac{1}{2}+\frac{k-1}{2m}}$$
and by Lemma 2.3,
$$\|z^{m-1}\oz^k \phi
u\|_{\frac{1}{2}+\frac{k-1}{2m}}\lesssim 
\|\oz^k \phi L u\|_\frac{k-1}{2m}+\|\oz^k \phi
u\|_{\frac{k-1}{2m}}+\|z^{k+2m-1}\phi'
u\|_{\frac{k-1}{2m}}.
$$
Now
$$\|\oz^k \phi u\|_{\frac{k-1}{2m}}\leq l.c.\|\phi
u\|_{-\frac{1}{2m}} + s.c. \|z^{m-1}\oz^k \phi
u\|_{\frac{1}{2}+\frac{k-1}{2m}}
$$
= RJ plus a term which can be absorbed by the
previous left hand side and a direct application of
Lemma 2.4 yields 
\begin{lemma} 
$$\|z^{k+2m-1}\phi'u\|_{\frac{k-1}{2m}}\lesssim
\|z^{2m}\phi' u\|_0+
\|z^{2k+m-2}\phi'u\|_{\frac{k-1}{m}-\frac{1}{2}}=RJ.$$
\end{lemma}

Putting these together,
$$\sum_{j=1}^{j=2}(LHS)_j \lesssim (LHS)_3+(LHS)_4 +
RJ.$$

\subsection{Estimation of $(LHS)_3$ and $(LHS)_4$}

Using the fact that $|F_z|\sim |z|^{2m-1},$ we have
$$(LHS)_3+(LHS)_4\equiv\|\phi\ol
u\|^2_{\frac{k-1}{2m}}+\|\phi\overline{z}^kLu\|^2_{\frac{k-1}{2m}}
$$
$$ = (\phi Pu, \phi u)_{\frac{k-1}{2m}} -([\overline{L},
\phi^2]|z|^{2k}Lu, u)_{\frac{k-1}{2m}}
- ([\phi^2,L]\overline{L}u, u)_{\frac{k-1}{2m}}$$
$$\leq l.c.\|\phi
P u\|^2_{\frac{k-1}{m}} + s.c. (LHS)_1
$$
$$+
2|(F_{\oz}\,\phi\phi'|z|^{2k}Lu, u)_{\frac{k-1}{2m}}| +
2|(F_z\phi\phi'\overline{L}u, u)_{\frac{k-1}{2m}}|
$$
$$\leq
l.c.\|\phi P u\|^2_{\frac{k-1}{m}} +
s.c. (LHS)_1 +
\{s.c.\|\phi
\oz^kLu\|^2_{\frac{k-1}{2m}}+l.c.\|z^{2m-1+k}\phi'
u\|_{\frac{k-1}{2m}}^2\}$$
$$ +
\{s.c.\|\phi z^{m-1}\ol
u\|_{\frac{k-1}{m}}^2+l.c.\|z^{m}\phi' u\|^2_0\}+RJ$$
$$\lesssim l.c.\|\phi P u\|^2_{\frac{k-1}{m}} + s.c.\{
(LHS)_1  +(LHS)_4+(LHS)_7\} +RJ$$
using Lemma 5.1 since the right hand side
(\ref{RJ1}) is $RJ$ for any $k.$
Thus, 
$$
\sum_1^4C_j(LHS)_j \leq 
s.c.(LHS)_7 + C_9\|\phi
Pu\|^2_{\frac{k-1}{m}} +  RJ.
$$
\subsection{Estimation of $(LHS)_5$ and $(LHS)_6$}
Setting $\sigma = \frac{k-1}{m}-\frac{1}{2}$ and using
Lemma 2.3,
$$(LHS)_5+(LHS)_6 \equiv \| z^{2k+m-2}\phi
u\|^2_{\frac{k-1}{m}} + \|\phi
z^{2k-1}Lu\|^2_\sigma$$
$$\lesssim
\|\phi z^{2k-1}Lu\|^2_{\sigma}+\|\phi 
z^{2k+m-3}u\|^2_{\sigma}+
\|z^{2k+3m-3}\phi'u\|^2_{\sigma}+RJ$$
$$\lesssim \|\phi z^{2k-1}Lu\|^2_{\sigma} + s.c.\|\phi
z^{2k+2m-3}u\|^2_{\frac{k-1}{m}}+\|\phi
u\|^2_{-\frac{2m-1}{2m}}+RJ
$$
$$\lesssim \|\phi z^{2k-1}Lu\|^2_{\sigma} +
s.c.(LHS)_5+RJ$$
(since $2k+2m-3\geq 2k+m-2$) with Lemma 2.4.

So we have to consider
$$(LHS)_6\equiv\|\phi z^{2k-1}Lu\|^2_{\sigma} \lesssim
|(\ol
\phi^2 |z|^{2(2k-1)} Lu,u)_\sigma|$$
$$+ \|\phi' z^{2(k+m-1)}u\|^2_\sigma + \|z^{2k-2}\phi
u\|^2_\sigma$$
$$\lesssim |(\ol
\phi^2 |z|^{2(2k-1)} Lu,u)_\sigma| +
s.c.\{(LHS)_1+(LHS)_5\} +RJ$$
by Lemma 2.4. 

Now for the inner product we have, again using Lemma
2.4,
$$|(\ol \phi^2 |z|^{2(2k-1)}Lu, u)_\sigma | \leq 
|(2\phi\phi'F_{\oz}|z|^{2(2k-1)}Lu, u)_\sigma |$$
$$+|(\phi \ol |z|^{2k}Lu, z^{2k-2}\phi u)_\sigma |+
|(\phi |z|^{4k-3}Lu, \phi u)_\sigma |$$
$$\lesssim s.c.\|z^{2k-1}Lu\|^2_\sigma +
l.c.\|z^{2k+2m-2}u\|^2_\sigma + \|\phi Pu\|^2_\sigma $$
$$+\|z^{2k-2}\phi u\|^2_\sigma + |(\phi L\ol u,
z^{2k-2}\phi u)_\sigma| + s.c.\|z^{2k-1}\phi Lu\|^2_\sigma
$$
so
$$|(\ol \phi^2 |z|^{2(2k-1)}Lu, u)_\sigma | \lesssim
s.c.\{(LHS)_6+(LHS)_5+(LHS)_1\}$$
$$+\|P u\|^2_\sigma+|(\phi L\ol u, z^{(2k-2)}\phi
u)_\sigma| + RJ$$
$$\leq s.c.\{(LHS)_6+(LHS)_5+(LHS)_1+(LHS)_8\}+\|P
u\|^2_\sigma +RJ.$$

Thus, so far, 
$$
\sum_1^6C_j(LHS)_j \leq 
s.c.\{(LHS)_7 + (LHS)_8 \}+ C_9\|\phi
Pu\|^2_{\frac{k-1}{m}} +  RJ
$$
but we still need to estimate both $(LHS)_7$ and $(LHS)_8$
since the definition of $RJ$ requires it. 

\subsection{Estimation of $(LHS)_7$ and $(LHS)_8$}

We proceed to estimate (a small multiple of) the following
expression $B,$ noting that $\|\phi L\ol u\|^2_\sigma
= \|L\phi \ol u\|^2_\sigma + RJ:$
$$B=\|z^{m-1}\phi \ol
u\|^2_{\frac{k-1}{m}}+\|\phi 
\ol\,\ol u\|^2_{\sigma}\;(+\|\phi 
L\ol u\|^2_{\sigma})\lesssim \|\phi L\ol u\|^2_\sigma +RJ $$
$$\lesssim |(\phi \ol L\ol u, \phi \ol
u)_{\sigma}| + |(\phi L\ol u,
\phi' F_z\ol u)_{\sigma}| +RJ
$$
$$= |(\phi \ol L\ol u, \phi \ol
u)_{\sigma}| + s.c.B+ RJ$$
$$\lesssim |(\phi \ol P u, \phi \ol
u)_{\sigma}| + |(\phi \ol\, \ol
|z|^{2k} L u, \phi \ol
u)_{\sigma}| + s.c.B + RJ$$
$$=B_1+B_2+s.c.B+RJ.$$
Now
$$ B_1=|(\phi \ol P u, \phi \ol
u)_{\sigma}|
\leq |(\phi P u, \phi L\ol
u)_{\sigma}|$$
$$ + |(F_{\oz}\, \phi\phi' P u, \ol
u)_{\sigma}|
\lesssim \|\phi P u\|^2_{\sigma} +
s.c.B + RJ
$$
while 
$$B_2=|(\phi \ol\, \ol |z|^{2k} L u, \phi\ol
u)_{\sigma}|\lesssim $$
$$\lesssim |(\phi \ol  |z|^{2k}\ol L u, \phi\ol
u)_{\sigma}| + |(\phi \ol z|z|^{2(k-1)}
L u, \phi\ol u)_{\sigma}|
$$
$$=B_{21}+B_{22}.$$
For $B_{21}$ we have 
$$B_{21}\lesssim |(\phi \ol  |z|^{2k}L\ol u, \phi\ol
u)_{\sigma}|+|(\phi \ol |z|^{2k}F_{z\oz} Tu,
\phi\ol u)_{\sigma}|$$
$$=B_{211}+B_{212}$$
with 
$$B_{211} \lesssim  
|(\phi  |z|^{2k}L\ol u, L\phi \ol
u)_{\sigma}|+|(F_{\oz}\,\phi' |z|^{2k}L\ol u,
\phi \ol u)_{\sigma}|
$$
$$\leq s.c.B + RJ.$$

For $B_{212},$ the bracket $[\,\ol, |z|^{2k}F_{z\oz}]$
will enter. This will contain a factor of
$z|z|^{2(k-1)+2(m-1)}$, and thus we may move
$\oz^{m-1}$ to the right hand side of the
inner product leaving a function 
$g(z)z|z|^{2(k-1)}z^{m-1}$ on the left, for a suitable
function $g.$  We have 
$$B_{212} = |(\phi \ol |z|^{2k}F_{z\oz} Tu,
\phi\ol u)_{\sigma}|$$
$$\lesssim |(\phi \ol F_{z\oz}|z|^{2k}  u, \phi \ol
u)_{\frac{k-1}{m}}| + |(\phi'\, \ol F_{z\oz}  |z|^{2k}
u, \phi \ol u)_{\sigma}|$$
$$\lesssim |(\phi F_{z\oz}|z|^{2k} \ol u, \phi \ol
u)_{\frac{k-1}{m}}| + |(\phi g(z)z|z|^{2(k-1)}z^{m-1} u,
\phi z^{m-1}\ol u)_{\frac{k-1}{m}}|$$
$$+ |(\phi F_{z\oz}|z|^{2k}  \ol u, \phi'\, \ol
u)_{\sigma}|+ |(\phi g(z)z|z|^{2(k-1)}z^{m-1}
u, \phi'  z^{m-1}\,\ol u)_{\sigma}|$$
$$\lesssim s.c. B + RJ +
\|\phi z^{2k+m-2}u\|^2_{\frac{k-1}{m}} \lesssim s.c.B +
(LHS)_5 + RJ 
$$

Finally, for $B_{22}$ we have 
$$B_{22}=|(\phi \ol z |z|^{2(k-1)} L u, \phi\ol
u)_{\sigma}|
$$
$$\lesssim |(\phi  z |z|^{2(k-1)} L u, L\phi\ol
u)_{\sigma}|+ |(\phi F_{\oz}z|z|^{2(k-1)} L
u, \phi'\ol u)_{\sigma}|$$
$$\lesssim \|\phi
z^{2k-1}Lu\|^2_{\sigma} + s.c.B + RJ = (LHS)_6+ s.c.B +
RJ.$$

Thus we have 
$$s.c.B \lesssim s.c.\{(LHS)_1+(LHS)_4+(LHS)_7+(LHS)_6\} +
\|\phi P u\|_{\frac{k-1}{m}}^2 + RJ$$
and hence 
$$\|\phi u\|_0^2 \lesssim \|\phi P u\|^2_{\frac{k-1}{m}}
+ RJ$$
or, iterating until $RJ$ is of arbitrarily low
order, for some $\tilde\phi \equiv 1$ near the support of
$\phi,$
$$\|\phi u\|_0^2 \lesssim \|\tilde\phi P
u\|^2_{\frac{k-1}{m}} + \|u\|^2_{-N}$$  

\section{Proof of Proposition 1.3 (Optimality)}
\renewcommand{\theequation}{\thesection.\arabic{equation}}
\setcounter{equation}{0}
\setcounter{theorem}{0}
\setcounter{proposition}{0}  
\setcounter{lemma}{0}
\setcounter{corollary}{0} 
\setcounter{definition}{0}
\setcounter{remark}{0}

\begin{proof}
For $\varphi$ of compact support, we set
$h_\lambda(z,t)=\varphi v_\lambda$ and 
$$v_\lambda
= exp(-\lambda(F-it-
(F-it)^2).$$ If $|z|$ is small enough, we have, from
above and below,
$$\Re \left(F-it-(F-it)^2 \right) = F-F^2+t^2
 \sim
|z|^{2m}+t^2,$$
and hence that $\|(\lambda (|z|^{2m}+t^2))^A
v_\lambda\|_\infty \sim C_A$ uniformly in $\lambda.$
Also, any function of compact support and equal to zero
in a neighborhood of the origin, such as a derivative
of a localizing function identically equal to one near
the origin, times $v_\lambda$ is of order
$\lambda^{-N}$ for any $N.$ 

Now $v_\lambda$ was chosen so that 
$\ol^Fv_\lambda=0,$ and, letting $H=F-it$ we
compute that 
$$L_{m,k}^F v_\lambda = -2\lambda
F_z(1+2H)v_\lambda$$
and hence that for some $A,$
$$\ol |z|^{2k} L v_\lambda = -2\lambda \ol(
|z|^{2k} F_z(1+2F-2it))v_\lambda\sim \lambda
|z|^{2k+2m-2}v_\lambda$$
$$\sim (\lambda
|z|^{2m})^\frac{2k+2m-2}{2m}
\lambda^{1-\frac{2k+2m-2}{2m}}v_\lambda \sim
\lambda^{-\frac{k-1}{m}}\left(\lambda
|z|^{2m}\right)^A.$$
Analogously, we have that as a principal term, 
$$\partial_t^s
v_\lambda\sim \lambda^s v_\lambda.
$$

Hence if there is an estimate of the form 
$$\|\psi v_\lambda\|_0 \lesssim \|\tilde \psi P_{k,m}^F 
v_\lambda\|_r + \|v_\lambda\|_{-\infty}$$
valid as $\lambda
\rightarrow \infty,$ for $\psi,\tilde\psi\in C_0^\infty,
\psi
\equiv 1$ near $0, \tilde\psi \equiv 1$ near
supp $\psi,$ then
$r\geq
\frac{k-1}{m},$ and an analogous argument holds for
Proposition 1.1. Finally the optimality at all levels
(other values of $s$) follows at once since the
vector field $\partial/\partial t$ commutes with the
differential operator $P_{k,m}^F.$
\end{proof}

\section{First proof of Theorem 1}
\renewcommand{\theequation}{\thesection.\arabic{equation}}
\setcounter{equation}{0}
\setcounter{theorem}{0}
\setcounter{proposition}{0}  
\setcounter{lemma}{0}
\setcounter{corollary}{0} 
\setcounter{definition}{0}
\setcounter{remark}{0}

As mentioned above, it is only in the `positive
cone' $\Gamma^+$ that hypoellipticity must be
shown, and then only in a small neighborhood of $z=0.$
Elsewhere the operator is subelliptic in the usual
sense or elliptic and (microlocal) hypoellipticity
with a gain of derivatives is well known. 

In $\Gamma^+,$ since $|\zeta|\lesssim |\tau|,$ showing
that high derivatives in $t$ exhibit the appropriate
gain will suffice. There is in general
no way (as yet) to pass from an estimate such as 
$$\|\phi \partial_t^r v\|_s \lesssim \|\tilde\phi
\partial_t^r Pv\|_{s+\frac{k-1}{m}} + \|v\|_{-\infty}$$
valid for {\it smooth} $v$ to the finiteness of $\|\phi
\partial_t^r u\|_s$ for a particular
(distribution) solution
$u$ to
$Pu=f$ with
$\tilde\phi\partial_t^rPu\in H^{s+\frac{k-1}{m}},$
however. In this situation, though, we may exploit the
fact that the coefficients of
$P$ do not depend on the variable
$t$ and introduce a cut-off function
$\chi(\tau)\equiv 1, |\tau| \leq 1$ and in
$C_0^\infty (|\tau|\leq 2), |\chi|\leq 1$ and set
$\chi_M(\tau) =
\chi(\tau/M).$ 

It is not difficult to see that with slight
modification, the above estimate may be applied  to
$v=\chi_M(\partial_t)
\gamma^+u$ where the operator $\chi_M(\partial_t)$ has
the obvious meaning (via Fourier transform) and
$\gamma^+$ is supported in $\Gamma^+$ and is equal to
one near the
$\tau$ axis. The slight modification is that we must
add a constant $C_u$ independent of $M$ to
handle derivatives of
$\gamma^+$. 

 Thus we
may write, suppressing $\gamma^+,$  
$$\|\phi\chi_M(\partial_t) \partial_t^r u\|_s
\lesssim
\|\tilde\phi\chi_M(\partial_t)
\partial_t^r Pu\|_{s+\frac{k-1}{m}} +
\| u\|_{-\infty} + C_u$$
and then let $M\rightarrow \infty$ to see that the
previous estimate holds also for the solution $u$ and
hence that $u$ is smooth in $t$ as well as in the other
variables. Strictly speaking, we would need to commute
$\phi$ (or $\tilde\phi$) with $\chi_M(\partial t),$
intruducing a term of lower order in $\partial t$
with a derivative on the localizing function, and hence
inductively handled. But this will not affect the
hypoellipticity. 

\section{Second proof of Theorem 1, 
$F(z,\oz)$ satisfying $(1.5)$}
\renewcommand{\theequation}{\thesection.\arabic{equation}}
\setcounter{equation}{0}
\setcounter{theorem}{0}
\setcounter{proposition}{0}  
\setcounter{lemma}{0}
\setcounter{corollary}{0} 
\setcounter{definition}{0}
\setcounter{remark}{0}

Here we will present a second proof of the
hypoellipticity of $P$ which extends
naturally to a proof of analytic hypoellipticity. It was
given in the model case $F(z,\oz) =
|z|^{2m}$ in \cite{Ta2005}. 

\subsection{The localization of powers of $\partial_t$}

Localization must be done very carefully, even
with $\phi$ depending on $t$ alone. For example, the
first 
 bracket $[L, \phi \partial_t^p]$ which we encounter
will contain 
$(L\phi (t)) \partial_t^p \sim iF_z\phi'\partial_t^p,$
which is problematic for any value of $m$ - no gain
in powers of $T$ and it is unclear even how to
estimate this expression. While
$\zeta$ may be small relative to $\tau,$ it is not
small relative to $|z|^{2m+1}\tau,$ so this error can
not be estimated easily.

In Kohn's work
\cite{K2005} ($m=1$), the 
$\oz$ (or $z$) in front of each derivative of $\phi$ is
carefully followed, and shown to provide, after
some work, a gain of $1/2$ derivative. Analyticity
was not considered in that paper, nor does is it
evident that it could be shown by those methods. 

Derridj and Tartakoff found in
\cite{DT2005} that an entirely different
approach, involving a delicately balanced localization
of
$\partial_t^p,$ led to analyticity rather directly, at
least for the case $F(z,\oz) = |z|^{2}.$ Then in
\cite{Ta2005} Tartakoff proved analyticity (and
$C^\infty$) hypoellipticity borrowing much of the
analysis of \cite{DT1988} for the case $F(z,\oz) =
|z|^{2m}.$ While replacing $2$ by $2m$ may not seem
like a big change, the degree of technical complexity
changed enormously. 

Here we look at $F(z,\oz)$ subject to the conditions
described above, namely
\begin{equation}F_{z\oz} = |z|^{2(m-1)}h,\, h(0) \neq 0,
F_z=\oz|z|^{2(m-1)}h \;{\hbox{ so that }} F_z =
\mu^{-1}\oz F_{z\oz}\end{equation}
in the language of \cite{DT1988} with $\mu$ real and
non-zero. 

To
introduce the general proof we include the short
argument from \cite{Ta2005}, namely the case of
$P_k, m=1,$ i.e., $F(z,\oz)=|z|^2.$

\subsection{The case of $P_k$ where $m=1$}

\begin{definition} 
For any pair of
non-negative integers, $(p_1,p_2),$ let
\begin{equation}\label{T} T = -i
\frac{\partial}{\partial t}\end{equation}
 and set
$$(T_0^{p_1,p_2})_\phi = \sum_{{a\leq p_1} 
\atop {b\leq p_2}}{\overline{L}^a\circ 
\overline{z}\,^a\circ 
T^{p_1-a}\circ
\phi^{(a+b)}\circ T^{p_2-b}\circ
{z}^b\circ L^b
\over a!b!},$$
where $$\phi^{(r)} =T^r\phi.$$
\end{definition}

\noindent Note that
the leading term (with
$a+b=0)$ is merely
$T^{p_1}\circ \phi \circ T^{p_2}$
which is equal to the operator 
$T^{p_1+p_2}$ on any open set $\Omega_0$ where
$\phi
\equiv 1.$

We have the precise commutation relations: 

\begin{proposition}\label{3.1} 
$$[L_0, (T_0^{p_1,p_2})_\phi] \equiv
(T_0^{p_1,p_2-1})_{\phi'}\circ L_0,$$
$$ [\overline{L_0},
(T_0^{p_1,p_2})_\phi] \equiv 
\overline{L_0}\circ (T_0^{p_1-1,p_2})_{\phi'},
$$
$$[(T_0^{p_1,p_2})_\phi,z] =
z\circ (T_0^{p_1,p_2-1})_{\phi'},
$$
and
$$
[(T_0^{p_1,p_2})_\phi,\overline{z}] =
(T_0^{p_1-1,p_2})_{\phi'}\circ \overline{z},
$$
\vskip.1in\noindent
where the $\equiv$ denotes modulo
$C^{p_1-p_1'+p_2-p_2'}$ terms of the form
\begin{equation}\label{error:|z|^2}
{\overline{L_0}^{p_1-p_1'}\circ
\overline{z}^{p_1-p_1'}\circ 
T^{p_1'}\circ\phi^{(p_1'+p_2'+1)}\circ
T^{p_2'}
\circ{z}^{p_2-p_2'}\circ
{L_0}^{p_2-p_2'}\over
(p_1-p_1')!(p_2-p_2')!}\end{equation}
 with either
$p_1'=0$ or
$p_2'=0,$ i.e., terms where all free $\partial_t$
derivatives have been eliminated on one side of
$\phi$ or the other. 
\end{proposition}
\begin{proof}
The proof is a straightforward calculation involving a
shift of index in the definition of
$(T_0^{p_1,p_2})_\phi.$
\end{proof}

The proof may be simple, but the result is remarkable:
whenever a bracket of the localization of
$(T^{p_1,p_2})_\phi$ with $L,\ol,z,$ or
$\oz$ occurs, the degree of $\partial_t$ drops by one,
the $L,\ol,z,$ or $\oz$ is not lost, and recursion is
always possible except when the `error' term
(\ref{error:|z|^2}) enters, where at least half of the
original $\partial_t$'s are converted into $L$ or $\ol$
- a decidedly favorable kind of error, and one which,
after additional use of the estimates, leads back to
localized powers of $T$ but of an order down
by a factor of $3/4$. 

\subsection{The general case}

We introduce the important vector field 
$$M=\mu^{-1}zL$$
has the property
$$[L, \overline{M}] \equiv F_z T \mod
\overline{M}$$ 
and hence 
$$[L, \phi \partial_t + \phi_t \overline{M}] \equiv 0
\mod (L\phi_t)\overline{M}.$$
In other words, we have managed to kill the most
disturbing term in the bracket of $L$ with $\phi
\partial_t.$ Note that addition of $\phi_t {M}$ would
not introduce new $\partial_t$ upon bracketing with $L$
and yet, upon bracketing with $\overline{L}$ would also
kill the most disturbing term. This suggests a
relatively straightforward generalization, along the
lines of \cite{DT1988} and \cite{BT2005}. 
\begin{definition} For $(p_1, p_2)$ as
 above, set 
$$\phi^{(d)} = 
\left(i\partial_t\right)^{d}\phi(t)$$
and set
$$N_b = \sum_{b'\leq b}
A^b_{b'}{M^{b'}\over b'!},$$ where the $A^b_{b'}$
(real) are to be determined subject to $A^b_b=1,$ and 
 $$\tilde{N}_a = \sum_{a'\leq a}
A^a_{a'}{{M^*}^{a'}\over
a'!}=(N_a)^*$$
where $M^*=-\ol \circ \oz \circ
\mu^{-1}$, and
set
$$(T^{p_1,p_2})_\phi = \sum_{{a\leq p_1} 
\atop {b\leq p_2}}\tilde{N}_a\circ  
T\,^{p_1-a}\circ
\phi^{(a+b)}\circ T\,^{p_2-b}\circ
N_b.$$

\end{definition}

We have
$$[\ol, (T^{p_1,p_2})_\phi] = [\ol, \sum_{a\leq
p_1,
b\leq p_2}\tilde{N}_a  T^{p_1-a}
\phi^{(a+b)} T^{p_2-b} N_b]$$
\begin{equation}\label{LTpphi}
=\sum_{a\leq
p_1,b\leq p_2}\bigg\{[\ol,
\tilde{N}_a] T^{p_1-a}
\phi^{(a+b)} T^{p_2-b} N_b
\end{equation}
$$-\tilde{N}_a T^{p_1-a} F_{\oz}
\phi^{(a+b+1)} T^{p_2-b} N_b
+ \tilde{N}_a  T^{p_1-a}
\phi^{(a+b)}
T^{p_2-b}[\ol, N_b]\bigg\}.$$

The last two terms on the right must cancel, to
preserve the balance, since both disturb the balance
between derivatives on $\phi$ and gain in powers of $T.$
We will choose the coefficients $A^b_{b'}$ of
${N}_b$ in such a way that, modulo
acceptable errors,   
\begin{equation}\label{recursion} [\ol,N_b] =
F_{\oz} T{N}_{b-1}.
\end{equation}
This will provide the needed cancellation via a shift
of index in $b$ in the sum just as in the case with
$F=|z|^2.$ The corresponding relation for brackets with
$L$ will follow by taking adjoints: again modulo
acceptable errors, 
\begin{equation}\label{recursion2}[L,\tilde{N}_a] = 
-\tilde{N}_{a -1}F_z T.
\end{equation}

Condition (\ref{recursion}) reads, using the
definition of ${N}_b,$ reads:
$$\sum_{b'=0}^b
A^b_{b'}{1\over b'!}\left[\ol, M^{b'}\right] =
-F_{\oz} T\sum_{b'=0}^{b -1} A^{b
-1}_{b'}{M^{b'}\over b'!}.$$

Expanding the brackets and keeping all factors of
$zL$ to the right, 
$$ {1\over b'!}\left[\ol, M^{b'}\right]={1\over
b'!}\sum_{1\leq b''\leq b'}{b'\choose
b''}ad_{M}^{b''}(\ol)M^{b'-b''}$$
$$={1\over
b'!}\sum_{1\leq b''\leq b'}{b'\choose
b''}ad_{M}^{b''-1}(-F_{\oz}T)M^{b'-b''}
$$
$$= -(F_{\oz}T)\sum_{1\leq b''\leq
b'}\frac{1}{b''!}
\,\frac{M^{b'-b''}}{(b'-b'')!}$$
since $M F_{\oz} = F_{\oz}.$
The condition (\ref{recursion}) thus requires, renaming
$b'-b''$ as $\overline{b}$ on the right just above,
\begin{equation}\label{rrB}\sum_{b''=1}
^{b-\tilde b} A^b_{\tilde b
+ b''}\frac{1}{b''!}
=A^{b-1}_{\tilde{b}}.
\end{equation}

Fortunately, we have investigated these
equations in \cite{DT1988} and, citing a result in the
book by Hirzebruch \cite{H1966} have explicit
solutions $\tilde{A}^*_*,$ unique
under the conditions that $\tilde{A}^q_0 = (-1)^q,$
namely 
$$\tilde{A}^r_s=\left(\left(\frac{t}
{e^t-1}\right)^{r+1}\right)^{(r-s)}(0)/(r-s)! $$

In addition, we will also need good expressions for
the other brackets: we compute  
$$[L, N_b] = [L,\sum_{b'=0}^b
A^b_{b'}\frac{M^{b'}}{b'!}] 
= \sum_{{1\leq b''\leq b'}\atop {b'\leq b}}
A^b_{b'}\,\frac{1}{b''!}\,
\frac{M^{b'-b''}}{(b'-b'')!}\circ
L$$
$$
[N_b, z] = z\circ \sum_{{1\leq b''\leq b'}\atop {b'\leq
b}} A^b_{b'}\,\frac{1}{b''!}\,
\frac{M^{b'-b''}}{(b'-b'')!}
$$
\begin{equation}\label{LoLN}
[\ol,
\tilde{N}_a] =
-\,\ol \circ \sum_{{1\leq a''\leq a'}\atop {a'\leq a}}
A^a_{a'}\,\frac{1}{a''!}\,
\frac{(M^*)^{a'-a''}}{(a'-a'')!},
\end{equation}
$$[\oz,\tilde{N}_a] =
\sum_{{1\leq a''\leq a'}\atop {a'\leq a}}
A^a_{a'}\,\frac{1}{a''!}\,
\frac{(M^*)^{a'-a''}}{(a'-a'')!}\circ \oz.
$$

In order to recognize these sums as $N$'s or
$\tilde{N}$'s, we need to be able to shift the lower
indices on $A^a_{a'}$
down by one. But this also we have done in
\cite{DT1988}, with the result that 
\begin{proposition} For any $r,s,$ and $c,$ 
$$A^r_s=\sum_{j=0}^{r-s}S_j^{r-s}\,A^{r-(c+j)}_{s-c}
$$
where
$$|S^k_\ell|\leq C^k.$$
\end{proposition}
These brackets, then, together with the Proposition,
immediately translate, setting $b''=c$ and
$\tilde{b}=b'-b'',$ into:
$$[L, N_b] = \sum_{\tilde{b}\leq b-c}
\frac{1}{c!}A^b_{\tilde{b}+c}\,\,
\frac{M^{\tilde{b}}}{\tilde{b}!}\circ
L=
\,
\sum_{\tilde{b}\leq
b-c-j}\frac{1}{c!}S^{c+j}_j\,A^{b-c-j}_{\tilde{b}}
\frac{M^{\tilde{b}}}{\tilde{b}!}\circ
L$$
or 
$$[L, N_b] = \sum_{{c+j\leq
b}\atop 1\leq c}\frac{1}{c!}S^{c+j}_j N_{b-c-j}\circ
L$$

Similarly, 
$$[\ol, \tilde{N}_a] = -\ol\circ \sum_{{c+j\leq
a}\atop 1\leq c}\frac{1}{c!}S^{c+j}_j
\tilde{N}_{a-c-j},$$ and 
$$[ \oz, \tilde{N}_a] = \sum_{{c+j\leq
a}\atop 1\leq c}\frac{1}{c!}S^{c+j}_j
\tilde{N}_{a-c-j}\circ \oz,$$ and
$$[{N}_b, z] = z\circ \sum_{{c+j\leq
b}\atop 1\leq c}\frac{1}{c!}S^{c+j}_j \tilde{N}_{b-c-j},$$

These precise commutation relations mean that the
whole localization $(T^{p_1,p_2})_\phi$ may be commuted
meaningfully with the vector fields $L, \ol$ and
with $z, \oz:$

\begin{proposition}\label{3.3} Modulo terms in which
either
$p_1$ or $p_2$ has been reduced to zero, and in view of
the cancellations ensured by (\ref{recursion}), 
\begin{equation}\label{LmTpphi}[L,
(T^{p_1,p_2})_\phi]\equiv
\sum_{{1\leq c, 0\leq j}\atop c+ j\leq
p_2}\frac{1}{c!}S^{c+j}_j
(T^{p_1,p_2-(c+j)})_{\phi^{(c+j)}}\circ
L\end{equation}
\end{proposition}
\begin{proof}
$$[L, (T^{p_1,p_2})_\phi]\equiv \sum_{{a\leq p_1} 
\atop {b\leq p_2}}\tilde{N}_a\circ  
T\,^{p_1-a}\circ
\phi^{(a+b)}\circ T\,^{p_2-b}[L, N_b]$$
$$\equiv \sum_{{a\leq p_1} 
\atop {b\leq p_2}}\sum_{{c+j\leq
b}\atop 1\leq c}\frac{1}{c!}S^{c+j}_j\tilde{N}_a\circ  
T\,^{p_1-a}\circ
\phi^{(a+b)}\circ T\,^{p_2-b} N_{b-c-j}\circ L$$
$$\equiv \sum_{{a\leq p_1} 
\atop {b\leq p_2}}\sum_{{c+j\leq
b}\atop 1\leq c}\frac{1}{c!}S^{c+j}_j\tilde{N}_a  
T\,^{p_1-a}
{\phi^{(c+j)}}^{(a+b-c-j)} T\,^{p_2-(c+j)-(b-c-j)}
N_{b-c-j} \circ L
$$
$$\equiv
\sum_{{1\leq c, 0\leq j}\atop c+ j\leq
p_2}\frac{1}{c!}S^{c+j}_j
(T^{p_1,p_2-(c+j)})_{\phi^{(c+j)}}\circ L$$
\end{proof}

Similarly we state, and omit the proofs, which are
virtually identical to that of the previous
proposition,  

\begin{proposition} 
\begin{equation}\label{LmbarTpphi} [\ol,
(T^{p_1,p_2})_\phi]\equiv -\ol\circ \sum_{{1\leq c,
0\leq j}\atop c+ j\leq p_1}\frac{1}{c!}S^{c+j}_j
(T^{p_1-(c+j),p_2})_{\phi^{(c+j)}}, 
\end{equation}
\begin{equation}\label{zbarTpphi}[\oz,
(T^{p_1,p_2})_\phi]\equiv
\sum_{{1\leq c, 0\leq j}\atop c+ j\leq
p_1}\frac{1}{c!}S^{c+j}_j
(T^{p_1-(c+j),p_1})_{\phi^{(c+j)}}\circ
\oz,\end{equation}
 and
\begin{equation}\label{zTpphi}[z,
(T^{p_1,p_2})_\phi]\equiv -\sum_{{1\leq c, 0\leq
j}\atop c+ j\leq p_2}\frac{1}{c!}S^{c+j}_j
(T^{p_1,p_2-(c+j)})_{\phi^{(c+j)}}\circ
z\end{equation}

\end{proposition}

What these commutation relations mean is that we may
move the vector fields of $P_m$ past
$(T^{p_1,p_2})_\phi$ freely, at each stage incurring
errors with the same vector fields and a gain in
derivatives in $(T^{p_1,p_2})_\phi.$ Thus we may
iterate the {\it a priori} inequality modulo errors of
nearly arbitrarily low order - all of the $\equiv$
signs above mean that we will ultimately arrive at
errors where either $p_1=0$ or $p_2=0.$  

 So we insert first 
$v=(T^{p_1,p_2})_\phi u$ into 
(\ref{ape}), then bring
$(T^{p_1,p_2})_\phi$ to the left of
$P=L\overline{L}+
 \overline{L}z^k\overline{z}^kL,$ and find that we have:
$$\|\overline{L}(T^{p_1,p_2})_\phi
u\|_0^2 +
\|\overline{z}^k{L}(T^{p_1,p_2})_\phi
u\|_0^2 +\|\Lambda^{-{k-1\over
2}}(T^{{p\over 2}, {p\over 2}})_\phi
u\|_0^2$$
\begin{equation}\label{est:ap2}\lesssim
|(P(T^{p_1,p_2})_\phi u, (T^{p_1,p_2})_\phi u)_{L^2}|
\end{equation}
$$\lesssim
|((T^{p_1,p_2})_\phi Pu, (T^{p_1,p_2})_\phi u)_{L^2}| + |([P,(T^{p_1,p_2})_\phi] u, (T^{p_1,p_2})_\phi u)_{L^2}|$$ and
by the above bracket relations, modulo the
same terms as above where all $T$'s from one side of
$\phi$ or the other have been `converted' into $L$'s or
$\overline{L}$'s, we have
$$([P,(T^{p_1,p_2})_\phi] u,
(T^{p_1,p_2})_\phi u) \equiv$$
$$= ([L\overline{L},(T^{p_1,p_2})_\phi]
u, (T^{p_1,p_2})_\phi u) + ([
\overline L z^k\overline{z}^k{L},(T^{p_1,p_2})_\phi] u, (T^{p_1,p_2})_\phi u)$$
$$= ([L,(T^{p_1,p_2})_\phi]\overline{L}
u, (T^{p_1,p_2})_\phi u) +
(L[\overline{L},(T^{p_1,p_2})_\phi] u,
(T^{p_1,p_2})_\phi u)$$
$$+([
\overline L ,(T^{p_1,p_2})_\phi]z^k\overline{z}^k{L} u, (T^{p_1,p_2})_\phi u)+ (\overline L[
 z^k,(T^{p_1,p_2})_\phi] \overline{z}^k{L}u, (T^{p_1,p_2})_\phi u)$$
$$+(\overline L z^k[
\overline{z}^k,(T^{p_1,p_2})_\phi]{L} u, (T^{p_1,p_2})_\phi u)+ 
(\overline L z^k\overline{z}^k[{L},(T^{p_1,p_2})_\phi] u, (T^{p_1,p_2})_\phi
u)$$
$$\equiv \sum_{{1\leq c, 0\leq j}\atop c+ j\leq
p_2}\frac{1}{c!}S^{c+j}_j(
(T^{p_1,p_2-(c+j)})_{\phi^{(c+j)}}
L\overline{L}u, (T^{p_1,p_2})_\phi u)$$
$$-\sum_{{1\leq c,
0\leq j}\atop c+ j\leq p_1}\frac{1}{c!}S^{c+j}_j(L\overline{L} 
(T^{p_1-(c+j),p_2})_{\phi^{(c+j)}} u, (T^{p_1,p_2})_\phi u)$$
$$- \sum_{{1\leq c,
0\leq j}\atop c+ j\leq p_1}\frac{1}{c!}S^{c+j}_j
(\overline{L}
(T^{p_1-(c+j),p_2})_{\phi^{(c+j)}} z^k\overline{z}^k
{L}u,  (T^{p_1,p_2})_\phi u)$$ 
\begin{equation}\label{est:ap3}
-\sum_{k'=1}^k \sum_{{1\leq c, 0\leq
j}\atop c+ j\leq p_2}\frac{1}{c!}S^{c+j}_j
(\overline L{z}^{k'}
(T^{p_1,p_2-(c+j)})_{\phi^{(c+j)}}
{z}^{k-k'}\overline z^k{L}u, 
(T^{p_1,p_2})_\phi u) \end{equation}
$$+\sum_{k'=0}^{k-1}\sum_{{1\leq c, 0\leq j}\atop c+ j\leq
p_1}\frac{1}{c!}S^{c+j}_j(\overline L{z}^k
\overline z^{k'} 
(T^{p_1-(c+j),p_1})_{\phi^{(c+j)}}
\overline z^{k-k'}{L}u, 
(T^{p_1,p_2})_\phi u)$$ 
$$+\, \sum_{{1\leq c, 0\leq j}\atop c+ j\leq
p_2}\frac{1}{c!}S^{c+j}_j(\overline L{z}^k
\overline z^k
(T^{p_1,p_2-(c+j)})_{\phi^{(c+j)}}{L}u,  (T^{p_1,p_2})_\phi u) 
$$ 
$$=A_1+A_2+A_3+\sum_{k'=1}^k A_{4,k'}+\sum
_{k'=0}^{k-1}A_{5,k'}+A_6.$$ 

Concerning the critical $L, \overline{L},z^k$ and
 $\overline{z}^k,$ note that in each term above
\begin{itemize}
\item no $L, \overline{L},$ power of $z$ or power
of $\overline{z}$ has been lost,
\item the order among $L, \overline{z}^k, z^k,$ and $
\overline{L}$ is preserved,
\item  letting $|q|=q_1+q_2,$ each term on the right
contains $(T^{q_1,q_2})_{\phi^{(|p|-|q|)}}$ with
$|q|<|p|$ and $|p|-|q|$ derivatives on $\phi,$
\item just as (\ref{est:ap3}) demonstrates the errors
which result in moving
$(T^{p_1,p_2})_\phi$ past the vector fields $L, \ol,
z^k\ol$ and $\oz^kL,$ further such brackets to position
the vector fields so as to make use the {\it a priori}
estimate again will produce similar errors, with
$|q|$ still lower and the `lost' $T$ derivatives
transferred to $\phi,$ 
\item iterating this process, together with a weighted
Schwarz inequality, will produce a sum of terms 
with $q_j\leq \frac{p}{2}$ of the
form 
$$l.c.\,\|\Lambda^\frac{k-1}{2m}
(T^{q_1,q_2})_{\phi^{(|p|-|q|)}}Pu\|_0^2
+s.c.\,\|\Lambda^{-\frac{k-1}{2m}}
(T^{\frac{p}{2},\frac{p}{2}})_\phi u)\|^2.$$
\item In fact, using larger constants
$\tilde S^{c+j}_j$ subject to the same kind
of bounds,
$|\tilde{S}^{c+j}_j|\leq
\tilde{C}^{c+j},$ we may replace all sums on the right
hand side above by suprema subject to the same range
restrictions on the indices.
\item This use of suprema allows us easily to iterate
everything on the right with easy control on the
constants until either $p_1$ or
$p_2,$ both of which start as $\frac{p}{2},$ drops to
zero, which may happen in two ways - either by
stepwise decrease as on the right hand side above from
successive brackets or by the single term in
Propositions \ref{3.1} and \ref{3.3} which is not
cancelled, the term with all $L$'s or $\ol$'s on one
side or the other in the definition of
$(T^{p_1,p_2})_\phi,$ whose principal term is 
$(T^{p_1,0})_{\phi^{(p_2)}} (zL)^{p_2}/p_2!$ or its
analogue with $p_1$ reduced to $0.$
\item At this point we no longer have an 
effective localization of powers of $T$ - for
example, brackets with $\ol$ are not corrected. We
proceed anyway, and when we lack a `good' vector field
such as $\ol$ (or of course $\oz^kL$), we create one by
integrating by parts:
\begin{equation}\label{lol}\|Lw\|^2 \lesssim \|\ol
w\|^2 + |(|z|^{2(m-1)}Tw,w)|\end{equation} 
to use up the $L$ and $\ol$ derivatives with the
byproduct of introducing up to half the 
number of new $T$ derivatives. 
\end{itemize}

Overall, then, the strategy has
been:
\begin{equation}\label{pto3p/2}
{\|\Lambda^{-\frac{k-1}{2m}}T^{p}u\|_{\{\phi
\equiv 1\}}}\rightarrow
{\|\Lambda^{-\frac{k-1}{2m}}(T^{\frac{p}{2},\frac{p}{2}})_\phi
u\|}\rightarrow\end{equation}
$$\rightarrow 
{\|\Lambda^{-\frac{k-1}{2m}}T^{\frac{3p}{4}}u\|_{supp
\phi}}\rightarrow 
{\|\Lambda^{-\frac{k-1}{2m}}T^{\frac{3p}{4}}u\|_{\{
\phi_1\equiv 1\}}}\rightarrow\ldots$$ for suitable
$\phi_1
\equiv 1$ on the support of $\phi.$ This will
continue, with a sequence of $\phi_j$ supported
in nested open sets as in \cite{Ta1978},
\cite{Ta1980} until only a negligible fraction
of $p$ is left, namely a bounded number of
derivatives. Since the order in $T^p$ is reduced by
a factor of $3/4$ each time, we will need
$\log_{4/3}p$ such nested open sets. Thus (where the
constant $C_{Pu}^{(p+\frac{k-1}{2m})}$ will
reflect bounds on the derivatives of $Pu$): 
\begin{equation}\label{pto3p/2}{\|\Lambda^{-\frac{k-1}{2m}}T^{p}u\|_{\{\phi
\equiv 1\}}}\
\lesssim
\|\Lambda^{-\frac{k-1}{2m}+3}u\|_{\{\phi_{\log_{{}_{4/3}}p}\equiv
1\}}+C_{Pu}^{(p+\frac{k-1}{2m})}\end{equation}
where the $3$ could be any other small integer. And
of course the whole derivation could have been done
at the
$H^s$ level: for any given $s,$ 
\begin{equation}\label{pto3p/2}
{\|\Lambda^{s-\frac{k-1}{2m}}T^{p}u\|_{\{\phi
\equiv 1\}}}\
\lesssim
\|\Lambda^{s-\frac{k-1}{2m}+3}u\|
_{\{\phi_{\log_{{}_{4/3}}p}\equiv
1\}}+C_{Pu}^{(p+\frac{k-1}{2m})}\end{equation}
which will end the story if this last norm is known
to be finite, provided  that the terms that arise
along the way are all similarly bounded. The most
important of these is of course 
$$|(\Lambda^{s+\frac{k-1}{2m}}(T^{\frac{p}{2},\frac{p}{2}})_\phi
Pu,
\Lambda^{s-\frac{k-1}{2m}}(T^{\frac{p}{2},\frac{p}{2}})_\phi
u)|$$ which shows that $Pu\in
H^{s+\frac{k-1}{2m}+p}$ in the largest of the
nested open sets implies that
$u\in H^{s-\frac{k-1}{2m}+p}$ in the smallest, a
loss of 
${k-1}$ derivatives.
 
The value of $s$ will be chosen so that
we know the norm on the right in (\ref{pto3p/2}) is
finite (for every distribution is locally in some
$H^{\tilde{s}}$), and then $p$ will be chosen so that
$Pu\in 
H^{s+\frac{k-1}{2m}+p}(\cup \phi_j)$
for that value of $s.$ It follows that $u\in
H^{s-\frac{k-1}{2m}+p}(\cap \phi_j).$

\begin{remark} For analyticity, one needs to ensure
that as we take $p$ larger and larger, the constants
satisfied by the Ehrenpreis-type localizers are
subject to bounds such that the estimate
(\ref{pto3p/2}) is uniform in $p.$ We have shown this
often before (cf. \cite{Ta1978}, \cite{Ta1980}) and the
arguments are the same here.
\end{remark}

\end{document}